\documentclass[10pt]{amsart}
\usepackage{amsmath,amsthm}

\newtheorem{theorem}{Theorem}[section]
\newtheorem{lemma}[theorem]{Lemma}
\newtheorem{proposition}[theorem]{Proposition}
\theoremstyle{remark}
\newtheorem*{remark}{Remark}

\newcommand{\R}{{\mathbb R}}

\newcommand{\supp}{\operatorname{supp}}
\newcommand{\intinf}{\int_{-\infty}^\infty}

\newcommand{\Nf}{N_f} 
\newcommand{\Nos}{S_f} 
\newcommand{\Zf}{Z_f} 

\renewcommand{\Re}{{\mathfrak{Re}}}
\renewcommand{\Im}{{\mathfrak{Im}}}

\renewcommand{\t}{\theta}

\renewcommand{\O}{{\mathcal{O}}}
\renewcommand{\i}{{\mathrm{i}}}
\renewcommand{\d}{{\mathrm{d}}}
\renewcommand{\^}{\widehat}

\newcommand{\E}{\mathbb{E}}
\newcommand{\aveTH}[1]{\left\langle#1\right\rangle_{T,H}} 
\newcommand{\LL}{\frac{\log T}{2\pi}}
\newcommand{\Om}{\Omega}
\newcommand{\weight}{w}

\begin{document}
\title{Linear statistics for zeros of Riemann's zeta function}
\author{C.P. Hughes and Z. Rudnick}
\address{Raymond and Beverly Sackler School of Mathematical Sciences,
Tel Aviv University, Tel Aviv 69978, Israel
({\tt hughes@post.tau.ac.il})}
\address{Raymond and Beverly Sackler School of Mathematical Sciences,
Tel Aviv University, Tel Aviv 69978, Israel
({\tt rudnick@post.tau.ac.il})}

\date{28 August 2002}

\thanks{Supported in part by  the EC TMR network 
``Mathematical aspects of Quantum Chaos'', EC-contract no
HPRN-CT-2000-00103.}

\begin{abstract}{
We consider a smooth counting function of the scaled zeros of the
Riemann zeta function, around height $T$. We show that the first few
moments tend to the Gaussian moments, with the exact number depending
on the statistic considered.
}
\end{abstract}

\maketitle


\section{Introduction}
In this paper we will examine linear statistics of zeros of the Riemann zeta function.  
Denote its nontrivial zeros by
$1/2+\i\gamma_j$, $j=\pm 1,\pm 2,\dots$ with $\gamma_{-j}=-\gamma_j$
and $\Re(\gamma_1)\leq \Re(\gamma_2)\leq \dots$.
Let $N(T)$ denote the number of zeros in the strip $0<\Re(\gamma)\leq T$, then $N(T)=\overline{N}(T)+S(T)$ where
\begin{align*}
\overline{N}(T) &= 1+\frac{1}{\pi}\Im\log\left(\pi^{-\i T/2}\Gamma(\tfrac{1}{4}+\tfrac{1}{2}\i T)\right)\nonumber\\
&= \frac{T}{2\pi}\log\frac{T}{2\pi e} +\frac{7}{8} +\O(1/T) 
\end{align*}

Selberg \cite{S2} has studied the distribution of the remainder term in the counting function, $S(t) = N(t)-\overline{N}(t)$, as $t$ varies between $T$ and $T+H$, where $H=T^a$ with $1/2<a\leq 1$. He showed that $S(t)$ has Gaussian moments in  the sense that when $T\to\infty$,
$$
\frac 1H \int_T^{T+H} 
\left| \frac{S(t)}{\sqrt{ (\log\log t)/2\pi^2}} \right|^{2k} \d t 
\to  \frac{(2k)!}{k!2^k}
$$

Fujii \cite{Fujii99}, among others, has studied the distribution of $N(t+h)-N(t)$ around $t$ near $T$. 
Since $\frac{1}{T}\int_0^T S(t)\;\d t\to 0$, the mean of this is asymptotic to $\overline N(T+h)-\overline N(T)$, and the error term, $S(t+h)-S(t)$, (which is thus asymptotically centered) has Gaussian moments, so long as $h$ is larger than the mean spacing of zeros at that height. That is, if $h\log T\to\infty$ but $h\ll 1$ then
\begin{equation*}
\frac{1}{T}\int_T^{2T} \left(\frac{S(t+h)-S(t)}{\sigma}\right)^{2k} \;\d t = \frac{(2k)!}{2^k k!} + \O\left(\frac{1}{\sigma}\right)
\end{equation*}
where
$$\sigma^2 = \frac{1}{\pi^2} \int_0^{h\log T} \frac{1-\cos t}{t}\;\d t$$
He has similar results for when $h\to\infty$ subject to $h\ll T$, when $\sigma^2$ is replaced by $\sigma^2=\frac{1}{\pi^2}(\log\log T -\log|\zeta(1+\i h)|)$.

Note that if $h$ is of the order of the mean spacing, that is if $h= \O( 1/\log T)$, then the main term is the same size as the error term (that is, both are $\O(1)$), and we may no longer conclude the distribution is Gaussian. This is not surprising, since for $h=\O(1/\log T)$ the distribution of $N(t+h)-N(t)$ in the large $T$ limit is discrete.

In this paper we will study the counting function in that critical scaling.
Rather than study $N(t)$ itself, instead we will investigate the distribution of a smooth version of the counting function in intervals of size comparable to the mean spacing, $2\pi/\log T$.
In particular, for a real-valued even function $f$, and real numbers $\tau$ and $T> 1$, set  
$$
\Nf(\tau): = \sum_{j=\pm 1,\pm, 2,\dots} 
f(\frac{\log T}{2\pi} (\gamma_j-\tau) )  \;.
$$
If $f$ is the characteristic function of an interval $[-1,1]$ and if all the $\gamma_j$ are real, then 
$\Nf(\tau)$ counts the number of zeros in the interval $[\tau-2\pi/\log T , \tau+2\pi/\log T]$.
However, we will take $f$ so that its Fourier transform,
$\^f(u):=\intinf f(x)e^{-2\pi\i xu}\;\d x$, is smooth and of compact support, and will not
assume the Riemann Hypothesis.  

As $T\to\infty$, we consider the fluctuations of $\Nf(\tau)$ as $\tau$ varies near $T$ 
in an interval of size about $H=T^a$, where $0<a\leq 1$.  
More precisely, given a weight
function $\weight\geq 0$, with $\intinf \weight(x)dx =1$, and $\^\weight$
compactly supported, we define an averaging operator 
$$
\aveTH{W}:= \intinf W(\tau)\weight(\frac{\tau-T}H)\frac{\d \tau}H \;.
$$

We will show that  the expected value of $\Nf$ is 
\begin{equation*}
\aveTH{\Nf} = \intinf f(x) \;\d x + \O(\frac{1}{\log T}) \; .
\end{equation*}

We will also show that for $\^f\in C_c^\infty(\R)$ the first
few moments of $\aveTH{(\Nf)^m}$ of $\Nf$ are Gaussian:
\begin{theorem}\label{thm:moments}
Let $H=T^a$ with $0<a\leq 1$, and let $\^f\in C_c^\infty(\R)$ be such that $\supp \^f \subseteq (-2a/m,2a/m)$. Then the first $m$ moments of $\Nf$ converge as
$T\to\infty$ to 
those of a Gaussian random variable with expectation $\intinf f(x) \;\d x$
and variance 
\begin{equation}\label{sigma_f}
\sigma_f^2 = \intinf \min(|u|,1) \^f(u)^2 \;\d u \;.
\end{equation}
\end{theorem}

The local statistics of the critically scaled zeros of the Riemann zeta function around height $T$ (that is, zeros scaled by the mean density, $\LL$) are believed \cite{Mont73,Odlyz} to behave like eigenangles of a random unitary matrix, when scaled by  $N/2\pi$, which is their mean density. Indeed, a similar result to the theorem above holds in random matrix theory \cite{hr1}: 
Since the $\t_n$ are angles, we consider the $2\pi$--periodic function
\begin{equation*}
F_N(\t) := \sum_{j=-\infty}^\infty f\left(\tfrac{N}{2\pi} (\t+2\pi j)\right)
\end{equation*}
and model $\Nf$ by
\begin{equation*}
\Zf(U) := \sum_{j=1}^N F_N(\t_j)
\end{equation*}
where  $U$ is an $N\times N$ unitary matrix with eigenangles $\t_1,\dots,\t_N$.

Writing $\E$ to denote the average over the unitary group with
Haar measure, then without any restrictions on the support of the
function $f$, we found in \cite{hr1} that $\E\{\Zf(U)\} = \int_{-\infty}^\infty f(x)\;\d x$,
and that the variance is
\begin{equation}\label{eq:rmt_var}
\sigma(f)^2 = \int_{-\infty}^\infty \min(1,|u|) \^f(u)^2\;\d u
\end{equation}
Observe that this is in complete agreement with the mean and
variance of $\Nf$ if $\^f$ has the same support
restrictions. 
Furthermore, we showed that for any integer $m\geq 2$, if 
$\supp \^f \subseteq [-2/m,2/m]$, then
\begin{equation*}
\lim_{N\to\infty} \E \left\{ \left( \Zf-\E\{\Zf\} \right)^m \right\} =
\begin{cases}
0 & \text{ if } m \text{ odd}\\
\frac{m!}{2^{m/2}  (m/2)!} \sigma(f)^m & \text{ if $m$ even}
\end{cases}
\end{equation*}
where $\sigma(f)^2$, the variance, is given in \eqref{eq:rmt_var}.

These are the moments of a normal random variable. However, the higher moments are not Gaussian, and we called this ``mock-Gaussian behaviour''

The random matrix results suggests that theorem~\ref{thm:moments} is not the complete truth. We expect the variance of $\Nf$, \eqref{sigma_f}, to hold without any restriction on the support of $\^f$, and the $m$-th moment of $\Nf$ to be Gaussian so long as $\supp \^f\ \subseteq [-2/m,2/m]$.

Random unitary matrices can also be used to model the low-lying zeros of Dirichlet $L$--functions, and mock-Gaussian behaviour was found there too \cite{hr1}. Other classical groups (symplectic, special orthogonal) are believed to model different classes of $L$--functions (like the quadratic $L$--functions), and they also show mock-Gaussian behaviour \cite{hr2}.

\section{Proofs}

\subsection{The explicit formula}

Set 
$\Om(r) = \frac{1}{2}\Psi(\frac{1}{4}+\frac{1}{2}\i r)+ \frac{1}{2}\Psi(\frac{1}{4}-\frac{1}{2}\i r) - \log\pi$, where $\Psi(s)=\frac{\Gamma'}{\Gamma}(s)$ is the polygamma function. We need the following version of Riemann's explicit formula:

\begin{lemma}\label{prop ef}
Let $g(u)\in C_c^\infty(R)$ and let $h(r)=\intinf g(u)e^{\i ru}\d u$. Then
\begin{equation}\label{explicit formula}
\sum h(\gamma_j) = h(-\frac \i 2)+h(\frac \i 2) + \frac 1{2\pi} \intinf h(r) \Om(r)\;\d r -\sum_{n=1}^\infty \frac{\Lambda(n)}{\sqrt{n}} \left( g(\log n)+g(-\log n) \right)
\end{equation}
where $\Lambda(n)$ is the von Mangoldt function.
\end{lemma}

For $\^f\in C_c^\infty(\R)$, setting
$$
h(r) = f(\LL(r-\tau)),\qquad 
g(u) = \frac {e^{-\i \tau u}}{\log T} \^f(\frac{u}{\log T})  
$$
we have $\Nf(\tau) = \overline{\Nf}(\tau) + \Nos(\tau)$ where 
\begin{multline}\label{mean W}
\overline{\Nf}(\tau)  = \frac 1{2\pi}\intinf f\left(\LL(r-\tau)\right) \Om(r) \d r \\
+  f\left(\LL(\frac {\i}{2}-\tau)\right) + f\left(\LL(-\frac {\i}{2}-\tau)\right)
\end{multline}
and 
\begin{equation}\label{\Nos}
\Nos(\tau) = -\frac 1{\log T} 
\sum_{n\geq 2}  \frac{\Lambda(n)}{\sqrt{n}} 
 \^f(\frac{\log n}{\log T})\left(e^{\i \tau \log n} + e^{-\i \tau \log n} \right)
\end{equation}

\begin{remark}
The conditions on $f$ (which are determined by the explicit formula), that its Fourier transform has compact support and is infinitely differentiable, can be considerably weakened to requiring that $f(r)$ is analytic in the strip $-c\leq \Im(r)\leq 1+c$ ($c>0$) such that $f(r) \ll (1+|r|)^{-(1+\delta)}$ ($\delta>0$, $r\in\R$).
\end{remark}

\subsection{The mean}

\begin{lemma}\label{lem:int f}
For all $A>1$,
\begin{multline}\label{asymp of mean{\Nf} }
\intinf f(\LL(r-\tau))\Om(r)\frac{\d r}{2\pi} \\
  = \frac{\Om(\tau)}{\log T}\intinf f(x)\d x  
+ \O( \frac 1{1+|\tau|} \frac 1{(\log T)^2} )+
\O\left(\frac{\log(1+|\tau|)}{(\log T)^A}\right)
\end{multline}
\end{lemma}

\begin{proof}
To evaluate the integral, we change variables and split the domain of
integration into two parts: 
\begin{align*}
\intinf f(\LL(r-\tau))\Om(r)\frac{\d r}{2\pi}  &= 
\frac 1{\log T} \intinf f(x)\Om(\tau+\frac {2\pi x}{\log T}) \d x \\ 
& = \frac 1{\log T} (\int_{|x|\leq Y} + \int_{|x|>Y})
f(x)\Om(\tau+\frac {2\pi x}{\log T}) \d x 
\end{align*}
where $Y\to\infty$ but $Y=o(\log T)$, say $Y=\sqrt{\log T}$. 

For the bulk of the integral $x/\log T$ is small, and so we expand 
$$
\Om(\tau+\frac{2\pi x}{\log T}) = \Om(\tau) + 
\O(\frac 1{1+|\tau|} \frac{|x|}{\log T})
$$ 
to find 
\begin{multline*}
\frac 1{\log T}\int_{|x|\leq Y} f(x)\Om(\tau+\frac {2\pi x}{\log T}) 
\d x \\
 =   \frac {\Om(\tau)}{\log T}\int_{|x|\leq Y} f(x)\d x 
+ \O\left(\frac 1{1+|\tau|} \frac 1{\log T}
\int_{|x|\leq Y}f(x)\frac{|x|}{\log T} \d x\right)  \\
 = \frac {\Om(\tau)} {\log T}\intinf f(x)\d x  + 
\O( \frac 1{1+|\tau|} \frac 1{(\log T)^2} )
\end{multline*}

For the tail of the integral we use $f(x)\ll |x|^{-N}$ for any 
$N\gg 1$ (which follows from $f\in C_c^\infty(\R)$) and
Stirling's formula (which yields $\Om(r) = \log(1+ |r|) +\O(1)$ for all $r\in\R$) to find that it is dominated by  
\begin{multline*}
\frac 1{\log T}\int_{|x|>Y} \frac 1{|x|^N} \log(1+|\tau|+\frac{|x|}{\log T})\; \d x \\
=\frac{2N}{Y^{N-1} \log T} \log(1+|\tau|+\frac{Y}{\log T}) + \frac{2N}{\log T} \int_{x>Y} \frac{1}{x^{N-1}} \frac{1}{(x+(1+|\tau|)\log T)}\;\d x\\
\ll \frac{\log(1+|\tau|)}{Y^{N-1}\log T }
\end{multline*}

Thus we have 
\begin{multline*}
\intinf f(\LL(r-\tau))\Om(r)\frac{\d r}{2\pi}
= \frac{\Om(\tau)} {\log T}\intinf f(x)\d x \\+ 
\O( \frac 1{1+|\tau|} \frac 1{(\log T)^2} ) + 
\O(\frac{\log(1+|\tau|)}{Y^{N-1}\log T})   
\end{multline*}
Taking $Y=\sqrt{\log T}$ gives \eqref{asymp of mean{\Nf} }. 
\end{proof}

\begin{lemma}\label{lem:mean}
For all  $\^f\in C_c^\infty(\R)$,
$$
\aveTH{\overline{\Nf}} = \intinf f(x)\d x  + \O(\frac 1{\log T}),\qquad T\to\infty
$$
\end{lemma}

\begin{proof}
Since $\Omega(\tau)=\log(1+|\tau|)+\O(1)$ for all $\tau$, we have
\begin{align*}
\aveTH{\frac{\Omega(\tau)}{\log T} \intinf f(x)\;\d x} &= \frac{\intinf f(x)\d x}{\log T}   
\left( \intinf \log(1+|\tau|) \weight( \frac{\tau-T}H ) \frac{\d \tau}{H}
+\O(1) \right)  \\
 &= \intinf f(x)\d x + \O(\frac 1{\log T} )
\end{align*}
and since the average of the error term in lemma~\ref{lem:int f} is similarly $\O(1/\log T)$, we have that
$$
\aveTH{\frac 1{2\pi}\intinf f\left(\LL(r-\tau)\right) \Om(r) \d r} = \intinf f(x)\d x + \O(\frac 1{\log T} )
$$

The averages of the polar terms, $ f(\LL(\frac {\i}{2}-\tau)) + f(\LL(-\frac {\i}{2}-\tau))$, in $\overline{\Nf}$ are bounded by $\O(\frac 1{H\log T})$ since by Parseval
$$
\intinf f(\frac{\log T}{2\pi}(\frac{\i}2 -\tau))\weight(\frac{\tau-T}H) \frac {\d
\tau}{H}  
 = \intinf \frac {2\pi}{\log T}\^f(-\frac {2\pi y}{\log T})e^{\pi y}
\^\weight(Hy)e^{-2\pi\i Ty}\d y
$$
and since $\^\weight$ has compact support, the integral is over $|y|\ll
1/H$ and is bounded by $\O(1/H\log T)$. 
\end{proof}

\begin{proposition}\label{prop:mean}
For $f$ with $\^f\in C_c^\infty(\R)$, if $H\to\infty$ then the mean value of $\Nf$ is given
by 
\begin{equation}\label{expected}
\aveTH{\Nf} = \intinf f(x)\d x +\O(\frac 1{\log T}), \qquad T\to\infty \;.
\end{equation}
\end{proposition}

\begin{proof}
In view of Lemma~\ref{lem:mean}, it suffices to show that the mean value of
$\Nos$ is zero as $H\to \infty$. Indeed, we have 
\begin{equation*}
\aveTH{\Nos} = \frac{-1}{\log T}
\sum_{n} \frac{\Lambda(n)}{\sqrt{n}}\^f(\frac{\log n}{\log T}) 
\left( \^\weight(\frac H{2\pi}  \log n) e^{-\i T \log n} 
+  \^\weight(-\frac H{2\pi}  \log n) e^{\i T \log n} \right)
\end{equation*}
Since $\^\weight$ has compact support, and the prime powers $n$ are at
least $2$, the summands vanish once $H$ is larger than a certain constant which depends upon the support of $\^\weight$.
\end{proof}

\subsection{The centered moments} 

\begin{proof}[Proof of Theorem \ref{thm:moments}.]
From lemma~\ref{lem:int f} and proposition~\ref{prop:mean} we have
\begin{align*}
 \aveTH{\left(\Nf-\aveTH{\Nf}\right)^m} &= \sum_{n=0}^m \binom{m}{n} \aveTH{\Nos^{m-n} (\overline{\Nf}-\aveTH{\Nf})^n}\\
&= \aveTH{\Nos^m}\left(1 + \O(\frac{1}{\log T})\right)
\end{align*}
and so it is sufficient to show that the $m$th moment of $\Nos$ is the
same as that of a centered normal random variable with variance given by \eqref{sigma_f}. This is achieved in the theorem~\ref{thm:moments of Nos}.
\end{proof}

Before we calculate the $m$-th moment of $\Nos$, as a warm-up we will consider the variance.
\begin{proposition}
If $H=T^a$ with $0<a\leq 1$ and $\supp \^f\subset (-a,a)$
then we have  
$$
\aveTH{(\Nos)^{2}} =  \intinf \min(|u|,1) \^f(u)^2 \d u +\O(\frac 1{\log T})
$$
\end{proposition}

\begin{proof}
Using the expression \eqref{\Nos}, multiplying out $(\Nos)^2$ 
and integrating we find 
\begin{multline*}
\aveTH{(\Nos)^2} =
\frac{1}{\log^2 T}
\sum_{\epsilon_1,\epsilon_2=\pm 1} \sum_{n_1,n_2}
\frac{\Lambda(n_1)}{\sqrt{n_1}} \frac{\Lambda(n_2)}{\sqrt{n_2}} \^f(\frac{\log n_1}{\log T}) \^f(\frac{\log n_2}{\log T})\\ 
\times \^\weight(\frac H{2\pi} (\epsilon_1 \log n_1 + \epsilon_2 \log n_2)) 
e^{-\i T (\epsilon_1 \log n_1 + \epsilon_2\log n_2)}
\end{multline*}

In order to get a nonzero contribution we need $\epsilon_1=-\epsilon_2$ (since once $H$ is larger than a certain constant which depends upon the support of $\^\weight$ all the $\epsilon_1=\epsilon_2$ terms vanish). Furthermore, since $\supp \^f\ \subseteq (-a,a)$ we have $n_1\leq T^{a-\epsilon}$ for some $\epsilon>0$, and therefore
$$
\^\weight(\frac {T^a}{2\pi} \log\frac{n_1}{n_2}) = 0
$$
(for large enough $T$) unless $n_1=n_2$.

Therefore, taking into account that
$\^\weight(0)=\intinf \weight(x)\d x = 1$ we find as soon as $T$ is sufficiently large,
$$
\aveTH{(\Nos)^{2}} = \frac{1}{\log^2 T}\sum_n 2\frac{\Lambda(n)^2}{n} \^f(\frac {\log n}{\log T})^2
$$
We note that by the Prime Number Theorem, as $T\to\infty$
\begin{align*}
\frac{1}{(\log T)^2}\sum_n 2\frac{\Lambda(n)^2}{n} \^f(\frac {\log n}{\log T})^2 & \sim 2\int_0^\infty u\^f(u)^2 \d u + \O(\frac 1{\log T})\\
&= \sigma_f^2 + \O(\frac 1{\log T})
\end{align*}
where $\sigma_f^2$ is given by \eqref{sigma_f}, this being true since $\supp\^f\subseteq (-1,1)$ by assumption of the theorem.
\end{proof}

\begin{theorem}\label{thm:moments of Nos}
If $H=T^a$ with $0<a\leq 1$ and $\supp \^f\subset (-\alpha,\alpha)$
then for 
$$
2\leq m<2a/\alpha
$$ 
we have  
$$
\aveTH{(\Nos)^{m}}\sim 
\begin{cases}
\frac{(2k)!}{k!2^k}\sigma_f^{2k},&m=2k\\0,&m=2k+1 
\end{cases} 
+\O(\frac 1{\log T})
$$
where the variance $\sigma_f^2$ is given by  \eqref{sigma_f}.
\end{theorem}

\begin{proof}
Using the expression \eqref{\Nos}, multiplying out $(\Nos)^m$ 
and integrating we find 
\begin{multline*}
\aveTH{(\Nos)^m} =
(\frac{-1}{\log T})^m
\sum_{\epsilon_1,\dots,\epsilon_m=\pm 1} 
\sum_{n_1,\dots,n_m} \prod_{j=1}^m
\frac{\Lambda(n_j)}{\sqrt{n_j}}\^f(\frac{\log n_j}{\log T})\\ 
\times \^\weight(\frac H{2\pi} \sum_{j=1}^m \epsilon_j \log n_j) 
e^{-\i T\sum_{j=1}^m \epsilon_j \log n_j}
\end{multline*}

Since $\^\weight$ has compact support, in order to get a nonzero
contribution  we need
$$
|\sum_{j=1}^m \epsilon_j \log n_j|\ll \frac 1H 
$$
Set $M=\prod_{\epsilon_j=+1} n_j$ and $N=\prod_{\epsilon_j=-1} n_j$. If $M\neq N$ then assume w.l.o.g. that $M>N$, say $M=N+u$ with $u\geq 1$. 
Thus for a non-zero contribution we need 
$$
\frac 1H\gg \log \frac{M}{N} = \log (1+\frac uN)\gg \frac 1N
$$
and hence $T^a = H\ll N\leq \sqrt{MN} \leq T^{m(\alpha-\epsilon)/2}$ since $n_j \leq T^{\alpha-\epsilon}$ by assumption on the support of $\^f$.
But $\alpha<2a/m$, so this is a contradiction. Therefore $M=N$, and $\sum  \epsilon_j \log n_j = 0$.

Thus for $T\gg 1$, we find (taking into account that
$\^\weight(0)=\intinf \weight(x)\d x = 1$)  
\begin{align}
\aveTH{(\Nos)^m} &= 
(\frac{-1}{\log T})^m \sum_{\epsilon_1,\dots,\epsilon_m=\pm 1} 
\sum_{\substack{n_1,\dots, n_m\geq 2\\ 
\sum_{j=1}^m \epsilon_j \log n_j = 0}}  
\prod_{j=1}^m\frac{\Lambda(n_j)}{\sqrt{n_j}}\^f(\frac{\log n_j}{\log T}) \nonumber\\
 &= \sum_{E \subseteq \{1,\dots,m\}} J(E) \label{eq:Sf in terms of J}
\end{align}
where 
\begin{equation}\label{eq:defn J}
J(E) = (\frac{-1}{\log T})^m 
\sum_{\substack{n_1,\dots, n_m\geq 2\\ \prod_{j\in E} n_j = \prod_{j\notin E} n_j}}  
\prod_{j=1}^m\frac{\Lambda(n_j)}{\sqrt{n_j}}\^f(\frac{\log n_j}{\log T})
\end{equation}
and the subset of indices $E$ denotes the $\epsilon_j$ which are positive. That is $j\in E$ iff $\epsilon_j=+1$.

Fix a subset $E\subset\{1,\dots,m\}$.
The sum in $J(E)$ is over tuples $(n_1,\dots,n_m)$ which satisfy $\prod_{j\in E} n_j = \prod_{i\notin E} n_i$.
We say that there is a {\em perfect matching} of terms if there  is
a bijection $\sigma$ of $E$ onto its complement $E^c$ in $\{1,\dots,m\}$ so that $n_j=n_{\sigma(j)}$, for all $j\in E$.
This can happen only if $m=2k$ is even and $\# E = \# E^c=k$.

Decompose
\begin{equation}\label{Decomposing J(E)}
J(E)=J_{diag}(E)+J_{non}(E)
\end{equation}
where $J_{diag}(E)$ is
the sum of matching terms - the diagonal part of the sum (nonexistent
for most $E$), and
$J_{non}(E)$ is the sum over the remaining, nonmatching, terms.

\subsubsection{Diagonal terms}
Assume that $m=2k$ is even.
There are $\binom{2k}{k}$
subsets $E\subset \{1,\dots 2k\}$ of cardinality $k=m/2$, and for each such subset $E$,
$J_{diag}(E)$ is the
sum over all $k!$ bijections $\sigma:E\to E^c$ of $E$ onto its complement,
of terms
$$
\left( \frac 1{(\log T)^2} \sum_n \frac
{\Lambda(n)^2}{n} \widehat f(\frac {\log n}{\log T})^2
\right)^k
$$
We evaluate each factor by using the Prime Number Theorem:
\begin{equation}\label{matching terms}
\begin{split}
\frac 1{(\log T)^2} \sum_n \frac
{\Lambda(n)^2}{n} \widehat f(\frac {\log n}{\log T})^2 
&\sim\frac 1{(\log T)^2} \int_2^\infty \frac{\log t}{t} \widehat f(\frac {\log t}{\log T})^2 \;\d t\\
&\sim\int_0^\infty u\widehat f(u)^2 \;\d u
\end{split}
\end{equation}
Since our function is even and supported inside
$(-\alpha,\alpha)$ and $\alpha<2a/m\leq 1$, we can  rewrite this as
$$
\frac 12 \int_{-\infty}^\infty \min(1,|u|) \widehat f(u)^2 \;\d u=:\sigma(f)^2/2
$$
This shows that for $m=2k$ even we have as $T\to \infty$  that
$$
\sum_{E\subset\{1,\dots ,m\}} J_{diag}(E) \to  \frac{(2k)!}{2^k k!} \sigma(f)^{2k }
$$

Below we will show that the nondiagonal terms $J_{non}(E)$ are
negligible, and hence by \eqref{eq:Sf in terms of J} and \eqref{Decomposing J(E)} we will
have thus proved Theorem~\ref{thm:moments of Nos}.
\end{proof}

\subsubsection{Bounding the off-diagonal terms $J_{non}(E)$}

We will show that
\begin{lemma}\label{bounding J_{non}}
$$
J_{non}(E) \ll \frac 1{\log T}
$$
\end{lemma}
\begin{proof}

Since
$$
\frac 1{\log T} \sum_{p}\sum_{ k\geq 3} \frac{\log p}{p^{k/2}} \ll
\frac 1{\log T} \sum_p \frac {\log p}{ p^{3/2}} \ll \frac 1{\log T}
$$
the contribution of cubes and higher prime powers to \eqref{eq:defn J} is negligible, and
we may assume in $J_{non}(E)$ that the $n_i$ are
either prime or squares of primes (upto a remainder of $\O(1/\log T)$). By the Fundamental Theorem of Arithmetic,
an equality $\prod_{j\in E}n_j = \prod_{i\in E^c}n_i$
forces some of the terms to match,
and unless there is a perfect matching of all terms,
the remaining integers satisfy equalities of the form $n_1n_2=n_3$
with $n_1=n_2=p$ prime and  $n_3=p^2$  a square
of that prime. Thus upto a remainder of $\O(1/\log T)$,
$J_{non}(T)$ is a sum of terms of the form
$$
\left( \frac 1{(\log T)^2}
\sum_{\substack{p\\k=1,2}} \frac {(\log p)^2}{p^k}
\widehat f(\frac{\log p^k}{\log T})^2 \right)^u \cdot
\left( \frac 1{(\log T)^3}
\sum_p \frac {(\log p)^3}{p^2} \widehat f(\frac{\log p}{\log T})^2
\widehat f(\frac{\log p^2}{\log T})
\right)^v
$$
with $2u+3v=m$, and $v\geq 1$.

We showed \eqref{matching terms} that the matching terms have an
asymptotic value, hence are bounded.
We bound the second type of term  by
$$
\frac 1{(\log T)^3}
\sum_p \frac {(\log p)^3}{p^2} \widehat f(\frac{\log p}{\log T})^2
\widehat f(\frac{\log p^2}{\log T}) \ll
\frac 1{(\log T)^3}
\sum_p \frac {(\log p)^3}{p^2} \ll \frac 1{(\log T)^3}
$$
Thus as long as $v\geq 1$  (that is if there is no perfect
matching of all terms), we get that the contribution is
$\O(1/\log T)$.
\end{proof}


%
\end{document}